\documentclass[12pt]{article}
\usepackage{mathrsfs}
\usepackage{amsfonts}
\usepackage{bbm}
\usepackage{amssymb}
\usepackage{amsmath,amssymb}
\newcommand{\al}{\alpha}
\newcommand{\be}{\beta}

\newcommand{\de}{\delta}

\newcommand{\mf}{\mathfrak}
\newcommand{\nin}{\noindent}

\newcommand{\ol}{\overline}

\newcommand{\seq}{\subseteq}

\newcommand{\si}{\sigma}

\newcommand{\vs}{\vspace*}
\newcommand{\vi}{\varphi}

\def\Lla{\Longleftarrow}
\def\Lra{\Longrightarrow}
\def \nin {\noindent}

\def \Lemma #1 {\vs{3mm}\nin {\bf Lemma #1} \it}
\def \Prop #1 {\vs{3mm}\nin {\bf Proposition #1} \it}
\def \Th #1 {\vs{3mm}\nin {\bf Theorem #1} \it}
\def \Cor #1 {\vs{3mm}\nin {\bf Corollary #1} \it}
\def \Ex #1 {\vs{3mm}\nin {\bf Example #1} \it}
\def \Proof {\vs{3mm}\nin {\bf Proof. }}
\def \part #1 {\hfil\break\hglue 12pt {\rm (#1)~}}

\def \qed {~\vrule height6pt width 6pt depth 0pt}
\def\fs{\footnotesize}

\title{
\bf\LARGE The structure of finite local principal ideal rings}
\author{{Tongsuo Wu$^a$\thanks{Corresponding author,\, tswu@sjtu.edu.cn (T. Wu)},  Houyi Yu$^b$\thanks{yhy178@163.com (H.Y. Yu)} and Dancheng Lu$^c$\thanks{ludancheng@suda.edu.cn (D. Lu)}}\\
  {\small $^{a,b}$Department of Mathematics, Shanghai Jiaotong University}\\
  {\small Shanghai 200240, P. R. China}\\
 {\small $^c$Department of Mathematics, Suzhou University, Suzhou 215006, P.R. China}\\
   \date{}
}
\begin{document}
\baselineskip=16pt
\maketitle

\begin{center}
\begin{minipage}{12cm}

 \vs{3mm}\nin{\small\bf Abstract.} {\fs A ring $R$ is called a PIR, if each ideal of $R$ is a principal ideal. A  local ring $(R,\mf{m)}$ is an artinian PIR if and only if its maximal ideal $\mf{m}$ is principal and has finite nilpotency index. In this paper, we determine the structure of a finite local PIR.}

\vs{3mm}\nin {\small Key Words:} {\small finite local ring, polynomial rings, ideal structure, ring extension }

\end{minipage}
\end{center}

\vs{3mm}\section{Introduction}

\vs{3mm}\nin   Throughout this paper, all rings $R$ are assumed to be commutative with identity $1_R$.  For a ring $R$,  let $U(R)$ be the set of invertible elements of $R$,   $J(R)$ the Jacobson radical of $R$ and $\mbox{char}(R)$ the characteristic of $R$. A ring $R$ is called a PIR if each ideal of $R$ is a principal ideal. For a local ring $(R,\mf{m})$ with a finite number of ideals (i.e., an artinian local ring $R$), $R$ is a PIR if and only if the maximal ideal $\mf{m}$ is cyclically generated (\cite[Proposition 8.8]{AM} or \cite[Theorem 2.1]{McL}). Furthermore, the number of nontrivial ideals of $R$ is $s$ if $s+1$ the nilpotency index of $\mf{m}$.   For any local ring $(R,\mf{m},\mf{k})$ with $\mbox{char}(\mf{k})=p$ and any pair $x,y\in \mf{m}$, note that $(x,y)=(x+my,x+ny)$ holds  for all integers $n,m$ satisfying $0\le m<n\le p-1$, thus $R$ has at least $p+2$ nontrivial ideals if $(x,y)$ is not cyclic: $(x,y),(y),(x+my)$ ($0\le m\le p-1$). Now for any integer $r\ge 4$, let $\mathbb I(r)$ be the set of finite local rings with $r$ nontrivial ideals. Then for any $(R, \mf{m},\mf{k})\in \mathbb I(r)$, if $\mbox{char}(\mf{k})\ge r-1$, then $R$ is a PIR. This implies that almost all finite local rings are PIR. In particular, a finite local ring $(R, \mf{m},\mf{k})$ with four nontrivial ideals is always a PIR if $\mbox{char}(\mf{k})\ge 3$.

Let $p$ be a prime number. Recall from \cite{McL}, page 255 (or \cite{Cohen}, Lemma 13, page 79 and Corollary 2, page 83.) that if $P$ is a field of characteristic $p$, then there exists
an unramified complete discrete valuation ring, $V$, {\it of characteristic $0$}
whose residue field is isomorphic to $P$. This ring $V$ is uniquely determined
up to isomorphism by $P$ and is called the {\it $v$-ring with residue field $P$}.  In \cite[Theorem 3.2]{McL}, the author uses the Cohen structure theory of complete local rings (see \cite{Cohen}) to prove that a local ring $(R,\mf{m})$ is an artinian PIR, with
$pR= \mf{m}^m$, and with $\mf{m}$ of nilpotency index $s$, if and only if
$R\cong V[X]/(f(X),X^s)$ for some Eisenstein polynomial $f(X)\in V[X]$ of
degree $m$, where $V$ is any $v$-ring with residue field $R/\mf{m}$.   Moreover, in \cite[Theorems 3.5, 3.6]{McL} the author classifies all tamely ramified (i.e., $p\nmid m$ in the aforementioned theorem) artinian local PIRs up to isomorphism by taking advantage of the established results on $v$-rings.

The purpose of this paper is to characterize a finite local PIR $(R,\mf{m})$ as some homomorphic image of $\mathbb Z_{p^r}[X,Y]$ and we determine the kernel, where $\mbox{char}(R)=p^r$. In Theorem 2.2,  for any finite local ring $R$, we construct a coefficient field of $R$ starting from a generator of the multiplicative cyclic group $(R/\mf{m})^*$. The main structural results are Corollary 2.3 and Theorem 3.1. Note that our structure theorems as well as the proofs are quite different from \cite{McL}.

\vs{3mm}\section{Finite local rings with a prime characteristic}

\vs{3mm} In this section, we classify all  finite local PIRs with a prime characteristic.  We start with a folk-law result:

\vs{3mm}\nin {\bf Lemma 2.1}    {\it If $(R,\frak{m},\mf{k})$ is an artinian local ring such that $\mf{k}$ is a finite field, then there exist a prime number $p$ and positive integers $r\le s\le t$ such that $char(R)=p^r$ and $|R|=p^t$, where $s$ is the  nilpotency index of $\mf{m}$.}

\Proof Assume that $\mf{m}$ has a nilpotency index $s$  and consider the following sequence of nonzero finite dimensional $\mathfrak{k}$-vector spaces
$$\mf{m}^{s-1},\mf{m}^{s-1}/\mf{m}^{s-2},\cdots, \mf{m}/\mf{m}^2.$$
If $char(\mf{k})=p$, then  $|R/\mf{m}|=p^m$ and  $|\mf{m}|=p^n$ for some natural numbers $m$ and $n$.  Clearly, $|R|=p^{m+n}$ and $n\ge m(s-1)$. Since $p^s=0$, it follows that $char(R)=p^r$ for some $r\le s$. \quad \qed

\vs{3mm} Now assume that $(R,\frak{m})$ is a finite  local PIR, where $\mbox{char}(R/\mf{m})=p$,  $\mf{m}=R\al$ and $\al$ has nilpotency index $s+1$. By Lemma 2.1, the characteristic of $(R,\mf{m})$ is $p^r$ for some $r$ with $1\le r\le {s+1}$. Then $\mathbb Z_{p^r}$ is a subring of $R$ and $p\in R\al$. Let $\varphi:\mathbb Z_{p^r}\to R/\mf{m}$ be the natural map and set $F=\mbox{im}(\varphi)$. Then $ F=\{\ol{0},\ol{1},\cdots, \ol{p-1}\}\cong \mathbb Z_p$, and $F$ is the prime subfield of the residue field $R/\mf{m}$. Assume further $R/\mf{m}= F[\ol{\be}]$ with $\be\in R$. Then
$$R=\mathbb Z_{p^r}[\be]+\mathbb Z_{p^r}[\be]\al+\mathbb Z_{p^r}[\be]\al^2+\cdots+\mathbb Z_{p^r}[\be]\al^{s}=\mathbb Z_{p^r}[\be,\al] $$
and there is a surjective ring homomorphism
$$\pi: \mathbb Z_{p^r}[X,Y]\to R=\mathbb Z_{p^r}[\be,\al], f(X,Y)\mapsto f(\be,\al).$$
Hence each finite local PIR is a proper homomorphic image of the polynomial ring  $\mathbb Z_{p^r}[X,Y]$, if $\mbox{char}(R)=p^r$. The problem is to determine the kernel of $\pi$. In this section, we consider the question for the case $\mbox{char}(R)=p$.

By \cite[Theorem 9]{Cohen} or \cite[Theorems 7.7, 7.8]{Eisenbud},  a complete noetherian local ring $(R, \mf{m})$ contains a coefficient field, if $R$ has the same characteristic as its residue field $R/\mf{m}$. It holds in particular for every finite local ring $R$ with characteristic $p$. For the ring $R$, we can construct a coefficient field of $R$ starting from any generator of the multiplicative cyclic group $(R/\mf{m})^*$.

\vs{3mm}\nin {\bf Theorem 2.2.}    {\it For a finite local ring $(R,\mf{m})$ with characteristic $p$, assume that $|R/\mf{m}|=p^m$.
Then there exist an element $\beta\in R$
and   a subring $A=\mathbb Z_p[\be]$ of $R$ such that $\beta^{p^{mt}-1}=1$ for some positive integer $t$ and $R/\mf{m}=\ol{A}$. In the case, $A$ is a  coefficient
field of $R$.  }

\Proof  Denote $R/\mf{m}$ by $F$ and set $k=p^m$. Then $F^*$ is a cyclic
multiplicative group of order $k-1$. Now let $\alpha$ be a generator of the group $F^*$. Then $\al^k=\al$, $\al^{k^t}=\al$ and hence $\al^{k^t-1}=1 $ holds for any positive integer $t$. Now let $s$ be the nilpotency index of $\mf{m}$ and choose $t$ such that $k^t\geq s$. Choose
$\beta_1\in R$ such that $\overline{\beta_1}=\alpha$.
Assume that $\beta_1^{k^t-1}=1+x$ for some $x\in \mf{m}$. Clearly
$$(1+x)^{k^t-1}=1+\sum_{r=1}^{s-1}C_{k^t-1}^rx^r.$$
Since $C_{k^t-1}^r\cdot x^r=C_{k^t}^r\cdot x^r+(-1)C_{k^t-1}^{r-1}\cdot x^r=(-1)C_{k^t-1}^{r-1}\cdot x^r$, it follows by induction that
$C_{k^t-1}^r\cdot x^r=(-1)^rx^r$. Then $(1+x)^{k^t-1}=1+\sum_{r=1}^{s-1}(-1)^rx^r $.
In this case, set $\beta=\beta_1(1+x)$. Then $\ol{\beta}=\al$ in $R/\mf{m}$ and $$\beta^{k^t-1}=\beta_1^{k^t-1}(1+x)^{k^t-1}=(1+x)(1+\sum_{r=1}^{s-1}(-1)^rx^r)=1.$$

Set $A=\{0,1,\beta, \beta^2, \cdots,\beta^{k^t-2}\}$. Clearly,
$R/\mf{m}=\ol{A}$.  We will show that
$A$ is a subring of $R$ and thus a field under the operations in $R$.
First, note that for any $\de\in R$, there
exist $\gamma\in A$ and $y\in \mf{m}$ such that $\de=\gamma+y$.
Consider the equation $X^{k^t}=X$. Clearly, every element in $A$ is a
root of the equation. If $\de\not\in A$, then there exist $\gamma\in
A$ and nonzero $y\in \mf{m}$ such that $\de=\gamma+y$, and so
$\de^{k^t}=(\gamma+y)^{k^t}=\gamma^{k^t}+y^{k^t}=\gamma=\de-y\not=\de$. This implies that $\de$ is not a
root of the equation, and  so for any $\de\in R$, $\de\in A$ if and
only if $\de$ is a root of the equation $X^{k^t}=X$. Given any $\de_1,\de_2\in A$, then
$(\de_1-\de_2)^{k^t}=\de_1^{k^t}+(-1)^{k^t}\de_2^{k^t}=\de_1-\de_2$, implying $\de_1-\de_2\in A$. Thus $A$ is a subring of
$R$ and $A=\mathbb Z_p[\be]$. Clearly, $A$ is also a field under the operation and thus $A=\{0,1,\beta, \beta^2, \cdots,\beta^{k-2}\}$. \quad \qed

\vs{3mm}  We now characterize finite local PIRs with characteristic $p$. The following corollary also follows from \cite[Theorem 3.1]{McL}, which is a direct consequence of Cohen's Theorem 9 (\cite{Cohen}).

\vs{3mm}\nin {\bf Corollary 2.3.}
 {\it Let $(R,\mf{m})$ be a finite local
 ring with } $\mbox{char}(R)=p$. {\it  Then $R$ is a PIR, $\mf{m}$ has nilpotency index $s$ if and only if $R\cong
F[X]/(X^s)$ for some
finite field $F$, where $X$ is an indeterminant over $F$. In particular, $R$ has two (respectively, three) nontrivial ideals if and only if $R\cong
F[X]/(X^s)$ for some
finite field $F$, where $s=3$ ($s=4$, respectively).}

\vs{3mm}\nin {\bf Proof.} Clearly,  $R$ is a PIR ring if $R\cong F[X]/(X^s)$ for any finite or infinite field $F$.

Conversely, assume that $R$ is a PIR and a local ring with characteristic $p$. Assume further $J(R)=R\al$, where $\al^s=0,\al^{s-1}\not=0$.  Let $A$ be the subfield
of $R$ as in Theorem 2.2. Then $R/R\al=\ol{A}$ and it implies $A+R\al^2+\cdots +R\al^{s-1}=R$. Thus  $R\al=A\al+A\al^2+\cdots +A\al^{s-1}$, and finally $R=A[\al]$.
Since $\al^s=0,\al^{s-1}\not=0$, it follows that $R\cong A[X]/(X^s)$.\quad $\Box$

\vs{4mm}\section{Finite local PIRs with characteristic $p^r$ $(r\ge 2)$}

\vs{3mm}Throughout this section, assume that $(R,\mf{m})$ is a local ring with $\mbox{Char}(R)=p^r$,  where $r\geq 2$ and $\mf{m}$ has nilpotency index ${s+1}$. Assume further $\mf{m}=R\al$, $R/\mf{m}=F[\ol{\be}]$ for some $\be\in R$, and $F\cong \mathbb Z_p$. Then $\al^{s+1}=0$ and $\al^s\not=0$.  We are going to determine the kernel of $\pi$, where
$$\pi: \mathbb Z_{p^r}[X,Y]\to R=\mathbb Z_{p^r}[\be,\al], f(X,Y)\mapsto f(\be,\al).$$

For any $f(X)=\sum_{i}a_iX^i$ in the polynomial ring $\mathbb Z_{p^r}[X]$, set $\ol{f}(X)=\sum_i\vi(a_i)X^i$, where $\varphi:\mathbb Z_{p^r}\to R/\mf{m}$ is the natural map.  Let $g(X)$ be a monic polynomial of $\mathbb Z_{p^r}[X]$ such that $\ol{g}(X)\in F[X]$
is the minimal polynomial of $\ol{\be}$ over $F$, where $\mbox{deg}(g(X))=\mbox{deg}(\ol{g}(X))$ and $F=\mbox{im}(\varphi)$.
For any monic polynomial $g_1(X)\in \mathbb Z_{p^r}[X]$, note that $\ol{g_1}(X)=\ol{g}(X)\in F[X]$ if and only if $g(X)-g_1(X)\in p\mathbb Z_{p^r}[X]$, which is denoted as $$g(X)\equiv g_1(X)(\mod p).$$
Since clearly there exists a natural ring epimorphism from $\mathbb Z_{p^r}[\be]$ to $F[\ol{\be}]=R/\mf{m}$ and the kernel is $p\mathbb Z_{p^r}[\be]+g(\be)\mathbb Z_{p^r}[\be]$, it follows that
$$\mf{m}=\mathbb Z_{p^r}[\be]p+Z_{p^r}[\be]g(\be)+\sum_{t=1}^s\mathbb Z_{p^r}[\be]\al^t,\quad U(R)=U(\mathbb Z_{p^r}[\be])+\mf{m}.\quad (\triangle)$$
These equalities will be used repeatedly in Section Four. Also, there exists a natural ring isomorphism $$\psi:\mathbb  Z_{p^r}[X]/(p,g(X))\to F[\ol{\be}]= R/\mbox{m},\, \ol{f(X)}\mapsto \ol{f}(\be).$$
In fact, consider $\si: \mathbb Z_{p^r}[X] \to F[\ol{\be}],\, f(X)\mapsto \ol{f}(\be).$ For any $f(X)$, write $$f(X)=g(X)h(X)+r(X)= g(X)h(X)+pr_1(X)+r_2(X),$$ where $\mbox{deg}(r_i(X))<\mbox{deg}(g(X)),$ and the coefficients of $r_2(X)$ is in $\{0,1,\cdots, p-1\}$. If $\ol{f}(\beta)=0$, then $\ol{r_2}(\be)=0$. Therefore $\ol{r_2}(X)=0$ in $F[X]$ and hence $r_2(X)=0$ holds in $\mathbb Z_{p^r}[X]$. This shows $\mbox{ker}(\si)=(p, g(X))$ and hence $\psi$ is an isomorphism .

 The following two claims  and an equality  will be used repeatedly in what follows:

\vs{3mm}\nin {\bf Claim $(i)$}: {\it For any $t(X)\in \mathbb Z_{p^r}[X]$,  if $t(X)$ is not in the ideal generated by $p$ and $g(X)$, then $t(\be)\in U(R)$.}

\vs{3mm}\nin{\bf Claim $(ii)$}: {\it $g(\be)\in \mf{m}$, hence $g(\be)^{s+1}=0$ and thus $\be$ is integral over $\mathbb Z_{p^r}$. }

$$\mathfrak m=\{0\}\cup U(R)\al\cup U(R)\al^2\cup\cdots\cup U(R)\al^s.\quad (*)$$

We now prove the main result of this section.

\vs{3mm}\nin {\bf Theorem 3.1.} {\it Assume that $p$ is a prime number and $r$ is an integer with $r\ge 2$. Then a ring $R$ is a finite local PIR with $\mbox{char}(R)=p^r$ and has exactly $s$ nontrivial ideals,  if and only if there exist integers $n\ge 1,m\ge 0$, a monic polynomial $g(X)$ in $\mathbb Z_{p^r}[X]$ and polynomials $u_i(X),v_j(X)\in \mathbb Z_{p^r}[X]\setminus (p,g(X))$ $(1\le i\le n,\, 1\le j\le m)$, where $\ol g(\ol{X})$ is irreducible in the polynomial ring $\mathbb Z_{p}[\ol{X}]$ (that is, in $\mathbb Z_{p^r}[X]/p\mathbb Z_{p^r}[X]$ ), $\mbox{deg}(g(X)=\mbox{deg}(\,\ol g(\ol{X})\,)$,\,  $\mbox{deg}(u_i(X))<\mbox{deg}(g(X))$ and $\mbox{deg}(v_j(X))<\mbox{deg}(g(X))$, such that
$R$ is isomorphic to  a factor ring $S/Q$, where $S=\mathbb Z_{p^r}[X,Y]$  and
$$Q=(pY^{s+1-t_1}, \, Y^{s+1},\, p-\sum_{i=1}^nu_i(X)Y^{t_i},\, g(X)-\sum_{j=1}^mv_j(X)Y^{s_j})$$
with $1\leq t_1< t_2< \cdots< t_n\leq s$, $1\leq s_1< s_2< \cdots< s_m\leq s$ and $(r-1)t_1\leq s<rt_1$.
}

\Proof
$\Lla$: Assume that $R= S/Q$. Clearly, $p^r\cdot R=0$, $Y^{s+1}\in Q$ and $Q\seq (p,g(X),Y)$. Also, there are natural ring isomorphisms
$$\mathbb Z_{p^r}[X,Y]/(p,g(X),Y)\cong \mathbb Z_{p^r}[X]/(p,g(X))\cong \mathbb Z_{p}[\ol{X}]/(\ol g(\ol{X})),$$
where $Z_{p}[\ol{X}]/(\ol g(\ol{X}))$ is a finite field by assumption on $g(X)$.
Set $\mf{m}=(p,g(X),Y)/Q$. Clearly, $R/\mf{m}\cong \mathbb Z_{p}[\ol{X}]/(\ol g(\ol{X})),\, \mf{m}=R\ol{Y}$, where $\ol{Y}^s\not=0$, $\ol{Y}^{s+1}=0$. It implies $p^{r-1}\cdot R\not=0$ whence $\mbox{char}(R)=p^r$. Then $\mf{m}$ is a nilpotent maximal ideal of $R$. Thus $R$ has a unique maximal ideal $R\ol{Y}$, which is principal and has nilpotency index $s+1$, and thus by Lemma 2.1, $R$ is a finite ring.

$\Lra:$ Let $(R,\mathfrak m,\mf{k})$ be a finite local PIR, where the nilpotency index is $s+1$ and $\mathfrak m=R\al$. Assume that $\mbox{char}(R)=p^r$, where $r\geq2$. Clearly, $p\in \mathfrak m$. By $(*)$ we  can suppose that $p=u\al^t$, where $u\in U(R)$ and $1\leq t\leq s$. Since $p^r\cdot 1_R=0$ while $p^{r-1}\cdot 1_R\neq 0$, it follow that $(r-1)t\leq s<rt$. We only need to determine the kernel of $\pi$, where
$$\pi: \mathbb Z_{p^r}[X,Y]\to R=\mathbb Z_{p^r}[\be,\al], f(X,Y)\mapsto f(\be,\al).$$

Consider the aforementioned ring isomorphism
$$\psi: \mathbb Z_{p^r}[X]/(p,g(X))\to F[\ol{\be}]= R/\mbox{m},\, \ol{f(X)}\mapsto \ol{f}(\ol \be).$$
There exists a polynomial $u_1(X)\in \mathbb Z_{p^r}[X]\setminus (p,g(X))$ such that $\psi(\ol {u_1(X)})=\ol u$. We can further assume that $deg(u_1(X))< deg(g(X))$. Then $u_1(\beta)-u\in \mathfrak m$, that is, $u=u_1(\beta)+m_1$ for some $m_1\in \mathfrak m$, hence  $p=u_1(\be)\al^{t}+m_1\al^{t}=u_1(\be)\al^{t_1}+u_2\al^{t_2}$ for some $u_2\in U(R)$ and some integer $t_2$, where $t_1=t$ and $t_2>t_1$.  For the same reason, there exists a  polynomial $u_2(X)\in \mathbb Z_{p^r}[X]\setminus (p, g(X))$ satisfying  $deg(u_2(X))< deg(g(X))$, such that $\psi(\ol {u_2(X)})=\ol u_2$. If further $u_2\al^{t_2}\neq 0$, then there exist elements $u_3\in U(R)$ and integer $t_3$ with $t_2<t_3$ such that $p=u_1(\be)\al^{t_1}+u_2(\be)\al^{t_2}+u_3\al^{t_3}$. Continuing the process, we obtain that
$$
p=u_1(\be)\al^{t_1}+u_2(\be)\al^{t_2}+\cdots+u_n(\be)\al^{t_n}
$$
where $u_i(x)\in \mathbb Z_{p^r}[X]\setminus (p, g(X))$ such that $deg(u_i(X))<deg(g(X)),\,\forall i=1,2,\cdots, n$, and $1\leq t=t_1<t_2<\dots<t_n\leq s$.

Next we consider the possible value of $g(\be)$. If $g(\be)\neq 0$, by $(*)$ we can assume $g(\be)=v\al^{s_1}$, where $1\leq s_1\leq s$ and $v\in U(R)$. Similarly, we have a polynomial $v_1(X)\in \mathbb Z_{p^r}[X]\setminus (p, g(X))$ such that $\psi(\ol {v_1(X)})=\ol v$ and $deg(v_1(X))< deg(g(X))$. Then $g(\be)=v_1(\be)\al^{s_1}+v_2\al^{s_2}$ with $v_2\in U(R)$ and $s_1<s_2$. If $v_2\al^{s_2}\neq 0$, then there exists a polynomial $v_2(X)\in \mathbb Z_{p^r}[X]\setminus (p, g(X))$ satisfying  $deg(v_2(X))< deg(g(X))$, such that  $\psi(\ol {v_2(X)})=\ol v_2$. Thus $g(\be)=v_1(\be)\al^{s_1}+v_2(\be)\al^{s_2}+v_3\al^{s_3}$ where $s_2<s_3$. Continuing the analysis, we obtain that
$$
g(\be)=v_1(\be)\al^{s_1}+v_2(\be)\al^{s_2}+\cdots+v_m(\be)\al^{s_m}
$$
where $v_j(X)\in \mathbb Z_{p^r}[X]\setminus (p, g(X))$ such that $deg(v_j(X))<deg(g(X)),\,\forall j=1,2,\cdots, m$ and $1\leq s_1<s_2<\dots<s_m\leq s$.

Take
$$Q=(pY^{s+1-t_1},\, Y^{s+1},\, p-\sum_{i=1}^nu_i(X)Y^{t_i},\, g(X)-\sum_{j=1}^mv_j(X)Y^{s_j})$$
Notice that if $g(\be)=0$, we substitute $g(X)-\sum_{j=1}^mv_j(X)Y^{s_j}$ by $g(X)$, and in this case, $m=0$.

Clearly, $Q\seq ker(\pi)$. Conversely, for any $f(X,Y)\in \mathbb Z_{p^r}[X,Y]$, write
$$
f(X,Y)\equiv Y^sw_s(X)+Y^{s-1} w_{s-1}(X)+\cdots+Y w_1(X)+w_0(X)(\mod Q) \quad (**)
$$
where $w_i(X)\in \mathbb Z_{p^r}[X]$ such that $deg(w_i(X))<deg(g(X))$.
Now assume that $f(\be,\al)=0$. Then $\al^sw_0(\be)=0$ and hence $w_0(\be)\in \mathfrak m$, that is, $w_0(X)\in (p, g(X))$, which together with $deg(w_i(X))<deg(g(X))$ yield that $w_0(X)=pw'_0(X)$. Thus we can substitute $w_0(X)$ by $pw_0(X)$ in $(**)$ and then we have
$$
f(X,Y)\equiv Y^sw_s(X)+Y^{s-1} w_{s-1}(X)+\cdots+Y w_1(X) (\mod Q) \quad (***)
$$
so that $\al^s w_1(\be)=0$, hence $w_1(x)=pw'_1(x)$ for some $w'_1(x)\in  \mathbb Z_{p^r}[X]$ with $deg(w'_1(x))<deg(g(x))$.
Thus we can substitute $w_1(X)$ by $pw_1(X)$ in $(***)$ and get that
$$
f(X,Y)\equiv Y^sw_s(X)+Y^{s-1} w_{s-1}(X)+\cdots+Y^2 w_2(X) (\mod Q).
$$
Since $s$ is finite, we can reach
$$
f(X,Y)\equiv 0(\mod Q).
$$
after finite steps. Therefore, $f(X,Y)\in Q$. This shows that $ker(\pi)=Q$, as required.\quad\qed

\vs{3mm} We remark that  by a similar discussion as above, we can get an alternative description for finite local
PIRs $R$ with $\mbox{char}(R)=p$, which is essentially the same with Corollary 2.3.

\vs{3mm}\nin {\bf Theorem 3.2.}
{\it Let $p$ be a prime number, and let $(R,\mf{m})$ be a finite local ring with $\mbox{char}(R)=p$. Then $R$ is a PIR with exactly $s$ nontrivial ideals if and only if there exist integers $m\ge 0$, a monic irreducible polynomial $g(X)$ in $\mathbb Z_{p}[X]$ and polynomials $v_j(X)\in \mathbb Z_{p}[X]\setminus (g(X))$ $(1\le j\le m)$, where  $\mbox{deg}(v_j(X))<\mbox{deg}(g(X))$, such that
$R$ is isomorphic to  a factor ring $S/Q$, where $S=\mathbb Z_{p}[X,Y]$  and
$$Q=(Y^{s+1},\, g(X)-\sum_{j=1}^mv_j(X)Y^{s_j})$$
with $1\leq s_1< s_2< \cdots< s_m\leq s$.}

\vs{3mm}\section{Finite local rings with at most three nontrivial ideals}

Note that in Theorem 3.1, the ideal $Q$ is generated by four polynomials. In some cases, it can be generated by less polynomials. In this section, we study the structure of finite local rings with at most three nontrivial ideals. It is easy to see that the maximal ideal of such a ring $R$ is principal and hence $R$ is a finite local PIR in each case. We will show that each $Q$ can be generated by at most three polynomials. The detailed verifications will be omitted in most cases.

Let $p$ be any prime number. For any positive integer $n$, it is well-known that there exists in $\mathbb Z_p[X]$ an irreducible monic polynomial $g(X)$ of degree $n$. In the following theorems, for any $f(X)=X^t+\sum_{i=0}^{t-1} a_iX^i\in \mathbb Z_p[X]$, set $f_r(X)= X^t+\sum_{i=0}^{t-1} a_iX^i\in \mathbb Z_{p^r}[X]$, where $1\le r\le 4$.

In \cite{Rein},  it is proved that a local ring $R$ has exactly one nontrivial ideal if and only if either $R\cong F[X]/(X^2)$ for some field $F$, or $R\cong V/p^2V$, where $V$ is a discrete valuation ring of {\it characteristic zero} and residue field of {\it characteristic} $p$, for some prime number $p$. For a finite local ring, use the discussions of Theorems 3.1 and 3.2 we can obtain the following result.

\vs{3mm}\nin {\bf Theorem 4.1.} {\it$R$ is a finite local ring  with exactly one nontrivial ideal and  residue field of characteristic $p$, if and only if $R$ is isomorphic to one of the  of following homomorphic images of a polynomial ring $\mathbb Z_{p^r}[X,Y]$ $(1\le r\le 2)$}:

$(1)$ {\it $\mathbb Z_{p}[X,Y]/(g(X), \,Y^2)$}

$(2)$ {\it $\mathbb Z_{p^2}[X]/(g(X)+pw(X))$}\\
\nin{\it where  $g(X)$  is monic and irreducible in $\mathbb Z_p[X]$ and $\mbox{deg}(w(X))<\mbox{deg}(g(X))$.}

\vs{3mm}\nin {\bf Theorem 4.2.} {\it $R$ is a finite local ring with exactly two nontrivial ideals and residue field of characteristic $p$,  if and only if $R$ is isomorphic to one of the following homomorphic images of a polynomial ring $\mathbb Z_{p^r}[X,Y]$ $(1\le r\le 3)$}:

$(1)$ {\it $\mathbb Z_{p}[X,Y]/(g(X), \,Y^3)$}

$(2)$ {\it $\mathbb Z_{p^2}[X,Y]/(g_2(X)+pw_2(X),\,p-u_2(X)Y^2,\, pY)$ or

 \hspace{0.6cm}$ \mathbb Z_{p^2}[X,Y]/(g_2(X)-v_2(X)Y+pw_2(X),\,p-u_2(X)Y^2,\,pY)$}

$(3)$ {\it $\mathbb Z_{p^3}[X]/(g_3(X)+pz_3(X)+p^2w_3(X))$
 }

\nin{\it where  $g(X)$ is a monic  and  irreducible polynomial of $\mathbb Z_p[X]$, $z(X),w(X),u(X)$ and $v(X)$ are polynomials of $\mathbb Z_p[X]$ such that $\mbox{deg}(x(X))<\mbox{deg}(g(X)),\forall x(X)\in \{z(X),w(X)\}$, $u(X)\not\equiv 0(\mod g(X))$ and $v(X)\not\equiv 0(\mod g(X))$.}

\vs{3mm} We give an outlined proof to the following:

\vs{3mm}\nin {\bf Theorem 4.3.} {\it $R$ is a finite local ring with exactly three nontrivial ideals and residue field of characteristic $p$,  if and only if $R$ is isomorphic to one of the following homomorphic images of a polynomial ring $\mathbb Z_{p^r}[X,Y]$ $(1\le r\le 4)$}:

$(1)$ $\mathbb Z_{p}[X,Y]/(g(X), \,Y^4)$

$(2)$ {\it $\mathbb Z_{p^2}[X,Y]/(g_2(X)+pw_2(X),\,Y^3-pu_2(X),\, pY)$ or

 \hspace{0.6cm}$ \mathbb Z_{p^2}[X,Y]/(g_2(X)-v_2(X)Y^2+pw_2(X),\,Y^3-pu_2(X),\,pY)$ or

 \hspace{0.6cm}$ \mathbb Z_{p^2}[X,Y]/(g_2(X)-v_2(X)Y+pw_2(X),\,Y^3-pu_2(X),\,pY)$ or

 \hspace{0.6cm}$ \mathbb Z_{p^2}[X,Y]/(g_2(X)-pYv_2(X),\,Y^2-pu_2(X)-pYz_2(X))$ or

 \hspace{0.6cm}$ \mathbb Z_{p^2}[X,Y]/(g_2(X)-pw_2(X)-pYy_2(X),\,Y^2-pu_2(X)-pYz_2(X))$ or

 \hspace{0.6cm}$\mathbb Z_{p^2}[X,Y]/(g_2(X)-v_2(X)Y-pYy_2(X)-pw_2(X),\,Y^2-pu_2(X)-pYz_2(X))$
}

$(3)$ {\it $\mathbb Z_{p^4}[X]/(g_4(X)+py_4(X)+p^2z_4(X)+ p^3w_4(X))$

}

\nin{\it where  $g(X)$ is a monic  and  irreducible polynomial of $\mathbb Z_p[X]$, $y(X), w(X),z(X),u(X)$ and $v(X)$ are polynomials of $\mathbb Z_p[X]$ such that $$\mbox{deg}(x(X))<\mbox{deg}(g(X)),\,\forall x(X)\in \{y(X),z(X),w(X)\},$$ $u(X)\not\equiv 0(\mod g(X))$ and $v(X)\not\equiv 0(\mod g(X))$.}

\Proof $\Lla:$ In  case (1) and (2) (respectively, case (3)), $\ol{(p,g(X), Y)}$ (respectively, $\ol{(p, g(X))}$) is the principal maximal ideal of $R$, and it has nilpotency index 4. Thus each ring is finite local and has three nontrivial ideals.

$\Lra:$ If $(R,R\al)$ is a finite local ring and $\al$ has nilpotency index four, then $\mbox{char}(R)\in \{p^i\,|\,1\le i\le 4\}$. Assume $R/R\al=\mathbb Z_p[\ol{\be}]$, then $R=\mathbb Z_{p^i}[\be,\al]$. Let $g(X)$ be the minimal polynomial of $\ol{\be}$ in $\mathbb Z_p[X]$. Then $g(\be)\in R\al$, where $$R\al=U(R)\al\cup U(R)\al^2\cup U(R)\al^3\cup\{0\}.$$ Consider the ring epimorphism
$$\pi: \mathbb Z_{p^i}[X,Y]\to R=\mathbb Z_{p^i}[\be,\al],\, f(X, Y)\mapsto f(\be,\al).$$

{\bf Case 1:} $\mbox{char}(R)=p$.

In this case, the result is stated in (1) and it follows from Corollary 2.3.

{\bf Case 2:} $\mbox{char}(R)=p^2$. In this case,  either $p\in U(R)\al^2$ or $p\in U(R)\al^3$.

{\it Subcase 2.1}: $p\in U(R)\al^3$. Under the assumption, we have $p\mf{m}=0$ and thus there exists a polynomial $u(X)\in \mathbb Z_p[X]\setminus (g(X))$ such that $Y^3-pu_2(X)\in \mbox{ker}(\pi)$.
If further $g_2(\be)\in U(R)\al^3\cup\{0\}$,  then
$$ \mbox{ker}(\pi)=(g_2(X)+pw_2(X),\,Y^3-pu_2(X),\, pY)$$
for some $w(X)\in \mathbb Z_p[X]$, where $\mbox{deg}(w(X))<\mbox{deg}(g(X))$. If further  $g_2(\be)\in U(R)\al^2$, then
$$ \mbox{ker}(\pi)=(g_2(X)-v_2(X)Y^2+pw_2(X),\,Y^3-pu_2(X),\,pY), $$
where $v(X)\not\equiv 0(\mod g(X))$ in $\mathbb Z_p[X]$. If $g_2(\be)\in U(R)\al$, then
$$ \mbox{ker}(\pi)=(g_2(X)-v_2(X)Y+pw_2(X),\,Y^3-pu_2(X),\,pY). $$

{\it Subcase 2.2}: $p\in U(R)\al^2$. In this case, it follows by $(\triangle)$ that $$Y^2-pu_{2}(X)-pYz_2(X)\in \mbox{ker}(\pi)$$ for some $u(X)\in \mathbb Z_p[X]\setminus (g(X))$, where either $z(X)=0$ or $z(X)\in \mathbb Z_p[X]\setminus (g(X))$.

If further $g_2(\be)\in U(R)\al^3$,  then by $(\triangle)$ we have  $g_2(X)-pYv_2(X)\in \mbox{ker}(\pi)$ for some $v(X)\in \mathbb Z_p[X]\setminus (g(X)) $, where $\mbox{deg}(v(X))<\mbox{deg}(g(X))$. Set
$$Q_4=(g_2(X)-pYv_2(X),\,Y^2-pu_2(X)-pYz_2(X)).$$
Then clearly $Q_4\seq \mbox{ker}(\pi)$,  $\{Y^4,pY^2,pg_2(X),Yg_2(X)\}\seq Q_4$ and $Y^3\equiv t_*(X)(\mod Q_4)$. Conversely, $\forall f(X,Y)\in \mbox{ker}(\pi)$,i.e., $f(\be,\al)=0$, we have
$$f(X,Y)\equiv Y^3t_*(X)+Y^2t_*(X)+Yt_*(X)+t_*(X)(\mod Q_4)$$
$$\equiv Yt(X)+s(X)(\mod Q_4)$$
 for some $t(X),s(X)\in \mathbb Z_{p^2}[X,Y],$ where $\mbox{deg}(t(X))<\mbox{deg}(g(X))$ and $\mbox{deg}(s(X))<\mbox{deg}(g(X))$. Then by Claim $(i)$, $\al t(\be)+s(\be)=0$ implies
$t(X)\equiv 0(\mod p)$ and $s(X)\equiv 0(\mod p)$. Hence $f(X,Y)\equiv p\cdot (Y\cdot x(X)+y(X))(\mod Q_4)$, where the degree of both $x(X)$ and $y(X)$ are less than $\mbox{deg}(g(X))$. It implies $Yy(X)\equiv 0(\mod Q_4)$, and therefore, $y(X)\equiv 0(\mod p)$. Then $f(X,Y)\equiv (pY)\cdot x(X)(\mod Q_4)$, and thus $x(X)\equiv 0(\mod p)$ by Claim $(i)$. Finally, $f(X,Y)\in Q_4$, and it implies $Q_4=\mbox{ker}(\pi)$ in the case.

If further $g_2(\be)\in U(R)\al^2\cup \{0\}$,  then it follows from $(\triangle)$ that
$$\mbox{ker}(\pi)=(g_2(X)-pw_2(X)-pYy(X),\,Y^2-pu_2(X)-pYz_2(X))$$ holds for some $y(X),z(X),w(X)\in \mathbb Z_p[X]$ with $$\mbox{deg}(x(X))<\mbox{deg}(g(X)),\,\forall x(X)\in \{y(X),z(X),w(X)\}.$$

If further $g_2(\be)\in U(R)\al$,  then it follows from $(\triangle)$ that
$$\mbox{ker}(\pi)=(g_2(X)-v_2(X)Y-pYy_2(X)-pw_2(X),\,Y^2-pu_2(X)-pYz_2(X)),$$
where $u(X),v(X),w(X),y(X)$ and $z(X)$ has the above mentioned properties.

{\bf Case 3:} $\mbox{char}(R)=p^3$.

The situation can not occur since $p\in U(R)\al\cup U(R)\al^2\cup U(R)\al^3.$

{\bf Case 4:} $\mbox{char}(R)=p^4$.

In this case, $p\in U(R)\al$, thus $\al$ is nothing new but a unital multiple  of $p$ in $\mathbb Z_{p^4}[X]$. In this case, it is easy to obtain $(3)$ by $(\triangle)$.\quad\qed

\vs{3mm}
In the following Propositions 4.4 and 4.5,  assume that $g_i(X)$  and $u_j(X)$  are monic polynomials in $\mathbb Z_{p^2}[X]$ such that $\mbox{deg}(u_r(X))<\mbox{deg}(g_1(X))$ ($r=1,3$), $\mbox{deg}(u_s(X))<\mbox{deg}(g_2(X))$ ($s=2,4$),
$g_i(\ol{X})$ are irreducible in $\mathbb Z_p[\ol{X}]$ (i.e., $\mathbb Z_{p^2}[X,Y]/(p, g_i(X))$) and $u_j(X)\not\equiv 0(\mod p)$.
We can assume further that the coefficients of $u_j(X)$ and $g_i(X)$ are in $\{0,1,\cdots,p-1\}$ ($1\le i\le 2,\,1\le j\le 4$). Then either $g_1(X)=g_2(X)$ or there exist $u(X),v(X)$ such that $$g_1(X)u(X)+g_2(X)v(X)\equiv 1(\mod p).$$

In Proposition 4.4,  set for $i=1,2$
$$
R_i=\mathbb Z_{p^2}[X,\,Y]/(g_i(X),\,Y^2-p\cdot u_i(X),\, pY),
$$
$$ \mf{m}_i=(p,\,g_i(X),\, Y)/(g_i(X),\,Y^2-p\cdot u_i(X),\, pY)$$
and $K_i=\mathbb Z_{p^2}[X]/(p,\,g_i(X))=\mathbb Z_p[\ol{X}]/(g_i(\ol{X}))$. Let $\pi_i: R_i\to K_i$ be the natural epimorphism. Note that $\mf{m}_i$ is the maximal ideal of $R_i$ for each $i$. Note also that $K_1\cong K_2$ if and only if $\mbox{deg}(g_1(X))=\mbox{deg}(g_2(X)).$

\vs{3mm}\nin {\bf Proposition 4.4.} (1) {\it If $R_1$ is isomorphic to $R_2$, then there exist a field isomorphism $\tau:K_1\to K_2$, a polynomial $v_2(X)$ in $\mathbb Z_{p^2}[X]\setminus (p,g_2(X))$ such that
$\ol{u_2(X)}=\ol{v_2(X)}^2\cdot \tau(\,\ol{u_1(X)})\,)$  holds in $K_2$.  }

(2) {\it Conversely,  if $g_1(X)=g_2(X)$ and there exist $v_i(X)$ in $\mathbb Z_{p^2}[X]$ such that $\ol{u_2(X)}=\ol{v_2(X)}^2\cdot\ol{u_1(X)}$ holds in $K_2$ and $\ol{u_1(X)}=\ol{v_1(X)}^2\cdot \ol{u_2(X)}$ holds in $K_1$  respectively, then $R_1$ is isomorphic to $R_2$. }

\Proof $(1)$ Assume that $\si: R_1 \to R_2 $ is a ring isomorphism.  Then $\si(\mf{m}_1)=\mf{m}_2$, and hence $\si$ induces a natural field isomorphism $\tau$ from $K_1$ to $K_2$. Clearly $\mf{m}_i=R_i\ol{Y}$, and then $\si(R_1\ol{Y})=R_2\ol{Y}$.
  Let $\ol{v(X,Y)}\cdot \si(\ol{Y})=\ol{Y}$ for some $v(X,Y)\in \mathbb Z_{p^2}[X,Y]$.  Note that $\ol{Y}^2=p\cdot\ol{u_i(X)}$ in $R_i$ ($i=1,2$). Thus in $R_2$
$$p\cdot \ol{u_2(X)}=\ol{v(X,Y)}^2\si(\ol{Y}^2)=p\cdot \ol{v(X,Y)}^2\si(\ol{u_1(X)}).$$
Note that  $p\cdot \ol{Y}=0$ in $R_2$, thus assume
$$p\cdot \si(\,\ol{u_1(X)}\,)= p\cdot\ol{ w(X)},\quad  p\cdot\ol{ v(X,Y)}=p\cdot\ol{ v_2(X)}  $$
for some $w(X),v_2(X)\in \mathbb Z_{p^2}[X]$. Then  $p\cdot[ u_2(\,\ol{X}\,)-v_2(\,\ol{X}\,)^2\cdot w(\,\ol{X}\,)]=0$ holds in the $cc$-local ring $(R_2,\,R_2\ol{Y})$, where $R_2=\mathbb Z_{p^2}[\ol{X},\ol{Y}]$.   By  $Claim\,(i)$, we must have
$u_2(X)-v_2(X)^2\cdot w(X)\in (p, \, g_2(X))$. Consider the following commutative diagram in which the isomorphism $\tau$ is induced by $\si$.

\[\begin{array}{ccc} R_1&\stackrel{\si}{\longrightarrow} & R_2 \\
\vcenter{\llap{ $\pi_1 $}}\,\Big\downarrow &  &\Big\downarrow\vcenter{\rlap{$\pi_2$}}\\
K_1 & \stackrel{\longrightarrow}{\tau} & K_2
\end{array}\]\par

\nin We get $\ol{u_2(X)}=\ol{v_2(X)}^2\cdot \tau(\,\ol{u_1(X)})\,)$ in $K_2$. Since $\ol{u_2(X)}\in U(R_2)$, both $\ol{v_2(X)}$ and $\tau(\,\ol{u_1(X)})$ are invertible in $R_2$. Note that $\ol{u_1(X)}=\ol{v_1(X)}^2\cdot \tau^{-1}(\,\ol{u_2(X)})\,)$ holds in $K_1$ for some $v_1(X)\in \mathbb Z_{p^2}[X]$.

$(2)$ By assumption, we can assume $\ol{w_1(X)}^2\cdot\ol{u_2(X)}=\ol{u_1(X)}$ in $K_2$.
Define a ring homomorphism $\vi_1: \mathbb Z_{p^2}[X,Y]\to R_2$ by $$1\mapsto 1,\,\vi(X)=X,\, \vi(Y)=\ol{w_1(X)Y} .$$
Since $p\cdot\ol{(p,g_2(X))}=0$ in $R_2$,
$$\vi_1(Y^2-pu_1(X))=\ol{w_1(X)Y}^2-p\ol{u_1(X)}=\ol{w_1(X)}^2\cdot (\ol{Y}^2-p\ol{u_2(X}))=0$$holds in $R_2$.
Then  $(g_1(X),\,Y^2-p\cdot u_1(X),\, pY)\seq \mbox{ker}(\vi_1)$, and $\vi_1$ induces a ring homomorphism $\ol{\vi_1}:R_1\to R_2$ such that
$$\ol{\vi_1}(\ol{X})=\ol{X},\quad \ol{\vi_1}(\ol{Y})=\ol{w_1(X)}\cdot \ol{Y}.$$
Under the assumption, there exists another ring homomorphism $\ol{\vi_2}:R_2\to R_1$ such that
$\ol{\vi_2}(\ol{X})=\ol{X},\quad \ol{\vi_1}(\ol{Y})=\ol{w_2(X)}\cdot \ol{Y}.$
Then $(\ol{\vi_2}\cdot \ol{\vi_1})(\ol{Y})=(\,\ol{w_1(X)}\cdot \ol{w_2(X)}\,)\cdot\ol{Y}$, and therefore $\ol{\vi_2}\cdot \ol{\vi_1}$ is a ring isomorphism.  \quad\qed

\vs{3mm}In the following Proposition 4.5, Set $R_1=S/\mf{I}_1,\, R_2=S/\mf{I}_2$, where $S=\mathbb Z_{p^2}[X,\,Y]$, $$\mf{I}_1=(g_1(X)-u_3(X)Y,\,Y^2-p\cdot u_1(X),pY)$$
$$\mf{I}_2=(g_2(X)-u_4(X)Y,\,Y^2-p\cdot u_2(X),\, pY).$$
Set $K_i=\mathbb Z_{p^2}[X]/(p,\,g_i(X))\cong \mathbb Z_p[X]/(\ol{g}_i(X))$. Let $\pi_i: R_i\to K_i$ be the natural epimorphism. Note that $K_1\cong K_2$ if and only if $\mbox{deg}(g_1(X))=\mbox{deg}(g_2(X)).$ By the proof of Proposition 3.2, we have

\vs{3mm}\nin {\bf Proposition 4.5.} (1) {\it If $R_1$ is isomorphic to $R_2$, then there exist a field isomorphism $\tau:K_1\to K_2$,  polynomials $v_2(X),\, u_5(X)$ in $\mathbb Z_{p^2}[X]\setminus (p,g(X))$ such that $\si(\,\ol{g_1(X)}\,)=\ol{u_5(X)}\cdot \ol{g_2(X)}\,$ holds in $R_2$ and $$\ol{u_2(X)}=\ol{v_2(X)}^2\cdot \tau(\,\ol{u_1(X)})\, $$
hold in $K_2$. }

(2) {\it Conversely,  if $g_1(X)=g_2(X)$ and there exist $v_i(X)$ in $\mathbb Z_{p^2}[X]$ such that $u_2(X)\equiv v_2(X)^2\cdot u_1(X)(\mod (p,g_2(X)^2)),\,u_4(X)\equiv v_2(X)u_3(X)(\mod (p,g_2(X)^2))$  and $u_1(X)=v_1(X)\cdot u_2(X)(\mod (p,g_1(X)^2)), u_3(X)\equiv v_1(X)u_4(X)(\mod (p,g_1(X)^2))$  respectively, then $R_1$ is isomorphic to $R_2$. }

\end{document}